\documentclass[11pt]{article}
\vspace{7cm} \textwidth 14cm \textheight 21.6cm
\usepackage{amsmath}
\usepackage{amsfonts}
\usepackage{url}
\usepackage[english]{babel}
\usepackage[dvipsnames]{xcolor}
\usepackage{hyperref}
\usepackage[pdftex]{graphicx}

\usepackage{changes}
\usepackage{latexsym,amsthm,amscd}
\newtheorem{thm}{Theorem}%[section]
\newtheorem{cor}{Corollary}
\newtheorem{prop}{Proposition}

\newtheorem{lem}{Lemma}
\newtheorem{Def}{Definition}

\newcounter{alphthm}
\theoremstyle{remark}
\usepackage{amssymb, graphicx, tcolorbox, wasysym,  pifont}
\numberwithin{equation}{section}

\newcommand{\RR}{\mathbb R}
\newcommand{\too}{\longrightarrow}
\newcommand{\be}{\begin{equation}}
	\newcommand{\ee}{\end{equation}}
\newcommand{\bc}{\begin{cor}}
	\newcommand{\ec}{\end{cor}}
\newcommand{\ben}{\begin{enumerate}}
	\newcommand{\een}{\end{enumerate}}
\newcommand{\beq}{\begin{eqnarray}}
	\newcommand{\eeq}{\end{eqnarray}}
\newcommand{\beqn}{\begin{eqnarray*}}
	\newcommand{\eeqn}{\end{eqnarray*}}

\newcommand{\bpf}{\begin{proof}}
	\newcommand{\epf}{\end{proof}}
\newcommand{\bl}{\begin{lem}}
	\newcommand{\el}{\end{lem}}
\newcommand{\bp}{\begin{prop}}
	\newcommand{\ep}{\end{prop}}
\newcommand{\bd}{\begin{Def}}
	\newcommand{\ed}{\end{Def}}
\newcommand{\bt}{\begin{thm}}
	\newcommand{\et}{\end{thm}}
\newcommand{\R}{I\!\! R} %$\R$

\title{Some  Rigidity Results on Complete Finsler
	Manifolds}

\author{ Asanjarani A.\footnote{E-mail: azam.asanjarani@auckland.ac.nz, ORCID ID:0000-0001-6115-073X} \,  ,\,\ Dehkordi H. R.\footnote{Corresponding author: E-mail: hengameh.r@ufabc.edu.br, Phone: +55(11) 49968332, ORCID ID:~0000-0002-1738-3373}\\
	\small $ ^* $ The University of Auckland, Auckland,  New Zealand\\
	\small $ ^\dagger $	Center of Mathematics, Computing and Cognition, Federal University of ABC,\\ Santo
	André, SP, Brazil}
\date{}

\begin{document}
	\maketitle
	
	\begin{abstract} 
		
		We provide an extension of Obata's theorem to Finsler geometry and establish some rigidity results based on a second order differential equation. Mainly, we prove that every complete connected Finsler manifold  of positive constant flag curvature is isometrically homeomorphic to an Euclidean sphere  endowed with a certain Finsler metric  and vice versa. Based on these results, we present a classification of  Finsler manifolds which admit a transnormal function. Specifically, we show that if a complete Finsler manifold admits a transnormal function with exactly two critical points, then it is homeomorphic to a sphere.

	\end{abstract}
	

	\section{Introduction}
	\label{intro}Rigidity describes  quite different concepts in
	mathematics. Historically, one of  the first rigidity
	theorems, proved by Cauchy in 1813, states that if the faces of a
	convex polyhedron were made of metal plates and the edges were
	replaced by hinges, the polyhedron would be rigid \cite{cauchy1813j}.
	Although rigidity problems were of immense interest to engineers,   the
	intensive mathematical study of these types of problems has occurred
	only in the late   20th century, see \cite{gruber1993handbook}. In geometry sometimes an
	object is considered as  rigid if it has flexibility and not elasticity. In other words,
	a geometrical rigidity implies invariant with respect to isometries. In Riemannian geometry,
	the sectional curvature is invariant under isometries. Hence, a space
	of positive constant curvature is transformed into the same space by
	each isometry. This fact is sometimes described as the ``strong
	rigidity'' of a space of constant curvature. 
	
	In  Finsler
	geometry, the encountered rigidity results are rather slightly
	weaker and they usually talk about under which assumptions on 
	the flag curvature -analogous to the sectional curvature in Riemannian
	geometry- the underlying Finsler structure is either Riemannian
	or locally Minkowskian. A famous treatise in this area is by
	Akbar-Zadeh \cite{akbar1988espaces} where he established 
	the following rigidity theorem for  compact manifolds: \emph{ Let $(M,g)
		$  be a compact
		without boundary Finsler manifold of constant flag curvature $K$.
		If $K<0$, then $(M,g)$ is Riemannian.
		If $K=0$, then $(M,g)$ is locally Minkowskian.}
	
	There are several papers in Finsler geometry with results similar to Akbar-zadeh's rigidity theorem but by considering different assumptions.
	Foulon addressed the case of strictly negative flag curvature  in
	Akbar-Zadeh's theorem.
	In \cite{foulon1997locally} he imposed  
	the additional hypothesis that the curvature is covariantly constant
	along a distinguished vector field on the homogeneous  bundle of tangent
	half lines to show that  the Finsler structure is Riemannian. 
	Also, he  presented a strong
	rigidity theorem for symmetric compact Finsler manifolds with
	negative curvature and proved that such manifolds are isometric to  locally
	symmetric negatively curved Riemannian spaces \cite{foulon2002curvature}. This
	extends  Akbar-Zadeh's rigidity theorem to a so-called ``strong
	rigidity'' one. 
	Shen \cite{shen2005finsler}   considered  the case of negative but not necessarily
	constant flag curvature by adding the assumption that
	the $S$-curvature is constant and showed that the Akbar-Zadeh's
	rigidity theorem still holds. 
	
	Following several rigidity theorems in the two joint
	papers \cite{kim2000rigidity} and \cite{kim2003finsler},   Kim  in \cite{kim2007locally} proved that:
	``Any compact locally symmetric Finsler manifold with positive
	constant flag curvature is Riemannian'. 
	Also,  Bidabad \cite{bidabad2007classification} established some rigidity theorems as an application of connection theory in Finsler geometry.
	Another rigidity result is presented by  Wu  \cite{wu2008some} who  proved that any locally symmetric Finsler manifold with nonzero flag curvature must be Riemannian.
	
	Finsler manifolds of positive flag curvature have been studied and classified by several researchers and a number of results have been generalized from Riemannian spaces of positive sectional curvature to Finsler manifolds of positive flag curvature, see for instance \cite{bao2002finsler, wilking2018revisiting, xu2015sp, xu2018delta}. Also,  Bidabad in \cite{bidabad2011compact}, using the same idea as in \cite{asanjarani2008classification}, provided a classification of simply connected compact Finsler manifolds. In  2018  Boonnam~et.~al. \cite{boonnam2018berwald} proved  that a complete \textit{Berwald manifold} with nowhere vanishing flag curvature must be Riemannian. Also, 
	several results and open problems about Finsler manifolds with positive curvatures are addressed  in \cite{deng2017recent}. 
	
	In this paper, we apply the \textit{adapted coordinate system} introduced in \cite{asanjarani2008classification} to study the strong rigidity of Finsler manifolds of positive constant flag curvature. Particularly, we show that: \emph{A complete $n$-dimensional  Finsler manifold  
		is  of positive constant flag curvature
		if and only if it is isometrically homeomorphic  to an $n$-sphere
		equipped with a certain Finsler metric.}
	This result complements the Akbar-Zadeh's rigidity theorem by considering the case of $K>0$. 
	
	Also, we provide an extension of Obata's theorem to Finsler geometry. Obata's theorem in Riemannian geometry says (see   \cite{wu2014note} for more details):  \textit{ Let $(M, g)$ be a complete connected  Riemannian manifold of dimension $n \geq 1$ which admits a non-constant smooth solution of Obata’s equation
		$\nabla dw + wg = 0$.
		Then $(M, g)$ is isometric to the $n$-dimensional round sphere $S^n$}.\\
	Here, we show that,\\
	{\textbf{Theorem~\ref{isometer2}}: \it{Let   $(M,g)$  be a complete connected Finsler manifold of
			dimension $n\geq2$.  In
			order that  there is  a non-trivial solution
			of
			$
			\nabla^{H}\nabla^{H}\rho+ C^2 \rho g=0
			$
			on $M$, it is necessary and
			sufficient that $(M,g)$ be isometric to an n-sphere of radius $1/C$.}}
	
	Further, we apply adapted coordinates to extend some results from the Riemannian transnormal functions to Finsler geometry. A Finsler transnormal  function is a natural generalization of distance functions. More precisely, a smooth function $\rho :M\to\mathbb{R}$ on a Finsler manifold  $(M,g)$ is called a Finsler transnormal  function if the Finsler norm of the gradient of $\rho$ is constant along each level set of $\rho$. 
	
	In the Riemannian geometry, transnormal functions have been studied for many years and some interesting results have been established, see for instance \cite{thorbergsson2000survey, wang1987isoparametric}. However, transnormal functions from the Finsler geometry point of view have received less attention. This is in spite of  several interesting problems that can  be tackled in this area and 
	applications of Finsler transnormal functions  in Physics, particularly in modeling the propagation of waves of wildfire and water, see   \cite{ dehkordi2020mathematical, dehkordi2019huygens, ekici2016note}. 
	To the best of our knowledge, the only works on Finsler transnormal functions are \cite{alexandrino2019finsler, he2016isoparametric}.  In \cite{alexandrino2019finsler} a generalization of some results of  \cite{wang1987isoparametric} to the Finsler geometry is presented
	and in \cite{he2016isoparametric} a classification of isoparametric functions on Randers-Minkowski spaces is presented.
	
	Here, we  extend the results of \cite{he2016isoparametric}
	to provide a classification of Finsler manifolds based on the number of critical points of a transnormal function defined on them:  \textit{If the transnormal function has no critical points, one critical point or two critical points then, respectively, it is conformal to the direct product of an open interval of the real line and some complete manifold, the Euclidean space, or the sphere.} Moreover, in Theorem~\ref{sph} we prove:  \textit{If the transnormal function on a compact Finsler manifold has exactly two critical points then the space is homeomorphic to the sphere.}

	The remainder of this paper is structured as follows. In Section \ref{pre}, we recall some basic definitions in Finsler geometry, including adapted coordinates for Finsler manifolds satisfying Eq.~\ref{c-field}, and, Finsler transnormal functions.
	In Section \ref{Sec:special}, we study a special case of Eq.~\ref{c-field} which is important  for establishing the main results of this work.  In section \ref{Otava}, we proceed with generalizing the Obata's theorem and  in Section \ref{trans}, we focus on Finsler transnormal functions and prove that any complete Finsler transnormal function with two critical points is homeomorphic to a sphere.

	\section{Preliminaries}\label{pre}
	In this section,  we review some definitions of Finsler geometry that we refer to through  this paper. More details can be found in \cite{shen2001lectures}. 
	\subsection{Finsler Manifolds}
	Let $M$ be a real n-dimensional  manifold of class $C^ \infty$ and $TM$ its tangent bundle, i.e. 
	$
	TM=\underset{x\in M}{{\bigcup}}\big\{(x,y): y\in T_pM\big\}.
	$
	A {\it{Finsler structure}} on $M$ is a function
	$F:TM \rightarrow [0,\infty )$, with the following properties: 
	\begin{itemize}
		\item [(i)] $F$ is smooth on the tangent bundle of non-zero vectors $TM_{0}:=TM\backslash \{0\}$;
		\item [(ii)]  $F$ is
		positively homogeneous of degree one in $y$, i.e.
		$F(x,\lambda y)=\lambda F(x,y),  \forall\lambda>0$, where $(x,y)$
		is an element of $TM$;
		\item[(iii)] The Hessian matrix of $F^{2}$, $(g_{ij}):=\left({1 \over 2}
		\left[ \frac{\partial^{2}}{\partial y^{i}\partial y^{j}} F^2
		\right]\right)$, is positive definite on $TM_{0}$.
	\end{itemize}
	A \textit{Finsler
		manifold} is a pair consisting of a differentiable manifold $M$ and
	a Finsler structure $F$ on $M$.  The tensor field $g$ with the components $g_{ij}$ is called the \textit{Finsler metric tensor} and we denote a Finsler manifold by $(M,g)$. 
	We denote the natural projection on $TM_0$ by  $ \pi$ and its differential  by $\pi_*$, i.e.  $\pi_*:TTM_{0}\rightarrow TM$. The \textit{vertical vector bundle} on $M$ is defined as $\displaystyle VTM:=\underset{ \nu \in  TM_0}{\bigcup} V_\nu TM$ where $V_\nu TM=\ker \pi_{*}$ is the set of vectors tangent to $\nu \in TM_0$. The complementary decomposition $HTM$ where  $HTM\oplus VTM=TTM_0$ is called the \textit{non-linear connection} on $TM_0$. The coefficients of the nonlinear connection are denoted by $G^i_j(x,y)$, where $\displaystyle G^i_j=\frac{\partial G^i}{\partial y^j}$ and $\displaystyle G^i=\frac{1}{4}g^{ik}(\frac{\partial^2F^2}{\partial y^k\partial x^j}y^j-\frac{\partial F^2}{\partial x^k})$.    By using the local coordinates $(x^{i},y^{i})$ on $TM$, called the \textit{line elements},  we have the local field of frames $\{\frac{\partial}{\partial x^{i}},\frac{\partial}{\partial y^{i}}\}$ on $TTM$. 
	Given a  non-linear connection, we can choose a  local field of frames  $\{\frac{\delta}{\delta x^i}\frac{\partial }{\partial y^i}\}$ on $TTM_0$ where $\frac{\delta}{\delta x^i}:=\frac{\partial}{\partial x^i}-G^j_i\frac{\partial}{\partial y^j}$
	and $\frac{\partial}{\partial y^{i}}$ are the set of
	vector fields on $HTM$
	and $VTM$, respectively.

	A 1-form of the \textit{Cartan connection} is given by $\displaystyle w^i_j=\Gamma^i_{jk}dx^k+C^i_{jk}\delta y^k$, where $\displaystyle \Gamma^i_{jk}=\frac12g^{ir}(\frac{\delta g_{rk}}{\delta x^j}+\frac{\delta g_{jr}}{\delta x^k}-\frac{\delta g_{jk}}{\delta x^r})$ and $\displaystyle C^i_{jk}=\frac12g^{ir}\frac{ \partial g_{jk}}{\partial y^r}$. The coefficients $\Gamma^i_{jk}$ and $C^i_{jk}$ are called coefficients of \textit{horizontal and vertical covariant derivatives} of the Cartan connection, respectively. Given a tensor field $T$ with the components $T^i_{jk}(x,y)$ on $TM$, the components of the Cartan horizontal covariant derivative of $T$, $\nabla^H T$, are given by 
	$$\nabla_rT^i_{jk}:=\frac{\delta}{\delta x^r}T^i_{jk}-T^i_{sk}\Gamma^s_{jr}-T^i_{js}\Gamma^s_{kr}+T^s_{jk}\Gamma^i_{sr}\,.$$

	Assume that  $\gamma:I\to M$ defined by $t\to x^i(t)$ be a smooth curve on $M$ and $\displaystyle \tilde{\gamma}(t)=\big (x^i(t),\frac{dx^i}{dt}\big)$ its natural lift on $TM$. We say that $\gamma$ is a \textit{geodesic} of the Finsler manifold $(M,g)$ if  $\nabla_{\dot{\tilde{\gamma}}} \dot{\tilde{\gamma}}=0$. Here, $\displaystyle\dot{\tilde{\gamma}}(t)=\frac{d}{dt}\frac{\delta y^i}{\delta x^i}+\frac{\delta y^i}{dt}\frac{\partial}{\partial y^i}$, where $\displaystyle \frac{\delta y^i}{dt}:=\frac{dy^i}{dt}+G^i_j\big(x(t),\frac{dx}{dt}\big)\frac{\partial x^j}{dt}$.\\
	\subsection{ Finsler Manifolds with a Non-trivial Solution  of $\nabla^{H}\nabla^{H}\rho= \phi g$ }
	Let $\rho:M\rightarrow \R$  be a scalar
	function on $M$ that satisfies  the following second order differential
	equation
	\be
	\label{c-field}
	\nabla^{H}\nabla^{H}\rho= \phi g,
	\ee
	where $\nabla^{H}$ is the Cartan horizontal covariant derivative and
	$\phi$ is a function of $x$ alone.  The connected component of a
	regular hypersurface defined by $\rho=constant$ is called a\emph{
		level set of $\rho$}. We denote by $\verb"grad" \rho$ the gradient
	vector field of $ \rho$ which is locally
	written in the form $\verb"grad" \rho =  \rho^{i}\frac{\partial}{\partial
		x^i}$, where  $\rho^i = g^{ij} \rho_j$, $\rho_j =\frac{\partial
		\rho}{\partial x^j}$ for $i, j, \dots \in \{1, \dots, n\}$.
	Note that the partial derivatives $\rho_j $ are defined on the manifold $M$
	while $ \rho^{i}$, the components of $\verb"grad" \rho$, are defined
	on its slit tangent bundle $TM_0$. Hence, $\verb"grad" \rho$
	can be considered as a section of $\pi^*TM\rightarrow TM_0$,
	the pulled-back tangent bundle over $TM_0$,
	and  its trajectories lie on $TM_0$.
	For more details see \cite{asanjarani2008classification} and references therein. 
	One can easily verify that the canonical projection of the trajectories of the vector field
	$\verb"grad" \rho $ are geodesic arcs on $M$ \cite{asanjarani2008classification}. Therefore, we can choose 
	local coordinates
	$(u^1=t, u^2,...,u^n)$ on $M$  such that $t$ is the parameter of the geodesic containing
	the  projection of a trajectory of the
	vector field $\verb"grad" \rho$ and the level sets of $\rho$ are given by $t=$constant. These geodesics are called \emph{
		$t$-geodesics}. Since in this   local coordinate system,
	the level sets of $\rho$ are given by $t=$constant,  $\rho$ may be considered as a function of $t$ only.
	In the sequel we will refer to these level sets and  these local coordinates as \emph{$t$-levels} and \emph{adapted
		coordinates}, respectively. Also, note  that along any $t$-geodesic, Eq.~(\ref{c-field}) reduces to the  second order differential equation 
	\be 
	\label{para}
	\frac{d^2 \rho}{d t^2}=\phi(\rho),
	\ee
	where $\phi$ is a function of $\rho$ which is differentiable at non-critical points.
	
	Let $(M,g)$ be a Finsler manifold and $\rho$  a
	non-trivial solution of Eq.~(\ref{c-field})  on $M$. Then, using  the  adapted coordinates, components of
	the Finsler metric tensor  $g$ are given by
	\begin{equation}
		\label{metrg}
		\footnotesize
		(g_{ij})= \left(\begin{array}{lcl}
			1 \quad {0\quad \ldots\qquad{ 0}}\\
			0 \quad {g_{22} \quad\ldots\quad g_{2n}}\\
			\vdots \qquad{ \quad \ldots}\\
			0 \quad {g_{n2} \quad\ldots \quad g_{nn}}
		\end{array}\right),
	\end{equation}
	and $t$ may be regarded as
	the arc-length parameter of $t$-geodesics.
	It can be easily verified that the Finsler metric form of $M$ is given by
	\be \label{meter}
	ds^{2}=(dt)^{2}+ \rho'^{2}f_{\gamma\beta}du^{\gamma}du^{\beta}, 
	\ee
	where  $f_{\gamma\beta}$ are components of a Finsler metric tensor
	on a $t$-level of $\rho$ and  $ \rho'^{2}f_{\gamma\beta}$ is the induced
	metric tensor of this $t$-level. Here,  prime denotes the
	ordinary differentiation with respect to $t$.
	In this
	paper, the Greek indices $\alpha, \beta, \gamma, \dots$ run over the
	range $2,3, \dots, n$.

	A point $o$ of $(M,g)$ is called a \emph{critical point} of $\rho$ if the
	vector field $\verb"grad" \rho$ vanishes at  $o$, or equivalently if $
	\rho'(o)=0$, see \cite{asanjarani2008classification}. 
	If a non-trivial solution of  Eq.~\eqref{c-field}  has some critical points, then $M$ possess some interesting properties. For instance:

	\setcounter{alphthm}{0}
	\begin{lem}\label{lem} \cite{asanjarani2008classification}
		Let $(M,g)$ be  an n-dimensional Finsler manifold which admits a non-trivial
		solution
		$\rho$  of Eq.~\eqref{c-field} with one critical
		point. Then any $t$-level set of $\rho$ with
		Finsler metric form $\overline{ds}^{2}= f_{\gamma\beta}du^{\gamma}
		du^{\beta}$, where $f_{\gamma\beta}$ is given by Eq.~(\ref{meter}), has the positive
		constant flag curvature  $\rho''^2(0)$.
	\end{lem}
	
	\begin{prop}\cite{asanjarani2008classification}\label{con}
		Let $(M,g)$ be a connected complete Finsler manifold of dimension $n\geq 2$. If $M$ admits a non-trivial solution of Eq.~(\ref{c-field}), then depending on the number of	critical points of $\rho$, i.e. zero, one or two respectively, it is conformal to
		\begin{itemize}
			\item [(a)] A direct product $J \times M$ of an open interval $J$ of the real line and an $(n-1)$-dimensional complete Finsler manifold $M$.
			\item[(b)]  An $n$-dimensional Euclidean space.
			\item[(c)] An $n$-dimensional unit sphere in an Euclidean space.
		\end{itemize}	
	\end{prop}
	\begin{lem}\cite{asanjarani2008classification}\label{sph2}
		Let $(M,g)$  be a simply connected and compact Finsler manifold of dimension $n>2$ which admits a solution of Eq.~(\ref{c-field}) with two critical points, then $M$ is homeomorphic to an $n$-sphere.
	\end{lem}
	
	%%%%%%%%%%%%%%%%
	\subsection{Transnormal Functions}\label{transs}

	Given a Finsler manifold $(M,g)$ and a  smooth function $\rho:M \too\RR$, if there exists a continuous function $\mathfrak{b}:\rho(M) \longrightarrow \RR$ such that 
	\begin{equation}\label{dis-tra}
		g(\verb"grad"  \rho,\verb"grad"  \rho)=\mathfrak{b}\circ \rho,
	\end{equation}
	then $\rho$ is called a {Finsler transnormal} function. It is not difficult to show that, given any vector $v$  tangent to $M$,  \begin{equation}\label{grad}
		g(v,\verb" grad" \rho)=d\rho v,
	\end{equation} see \cite{shen2001lectures} for details. Recall that  a critical point  is a point $o$  such that $\rho'(o)=0$, we define  a regular point as a point of $M$  which is not critical. The regular and critical values are  images of regular and critical points, respectively, under $\rho$. The connected component of the pre-image of a regular value, $\rho^{-1}(t)$, is called a\emph{
		regular	level set of $\rho$} and the connected component of the pre-image of a critical value is called a\emph{
		singular	level set of $\rho$}. From Eq.~(\ref{dis-tra}), one deduces that the function $\mathfrak{b}$ is smooth on $\rho(M^0)$, where $M^0$ is the subset of $M$ containing the regular points  \cite{dehkordi2018finsler}. 
	
Given a Finsler manifold $(M,g)$ and any two points $p, \,q\in M$, the \textit{Finsler distance} from $p$ to $q$  is  defined as
\begin{equation}\label{dis}
	d(p,q):=\inf_{{\gamma}}\int_{a}^{b}\sqrt{g({\gamma^\prime}(t), {\gamma^\prime}(t))}dt,
\end{equation}
where the infimum is taken over all piece-wise smooth curves ${\gamma}:[a,b] \longrightarrow M$ joining $p$ to $q$.	One special example of Finsler transnormal functions  is the Finsler distance function: Given a compact subset $A\subset M$,  the Finsler distance function from $A$ to $p$ is given by $\rho:M\too\RR$ where $\rho(p)=d(A,p)$. One can prove that $\rho$ is locally Lipschitz continuous \cite{shen2001lectures} and therefore it is  differentiable almost everywhere. Also, it is not difficult to show that  the Finsler distance function $\rho$ satisfies $g(\verb"grad"  \rho,\verb"grad"  \rho)=1$ (see  Lemma $3.2.3$ of \cite{shen2001lectures}).  So, the Finsler distance function associated to a given  Finsler manifold is a transnormal function with  $\mathfrak{b}=1$ in~ Eq.~(\ref{dis-tra}).

	Some interesting properties of Finsler transnormal functions are presented in \cite{alexandrino2019finsler}. For instance, 
	\begin{prop}\cite{alexandrino2019finsler}
		\label{proposition-equidistant}
		Let $(M,g)$ be a Finsler manifold and $\rho: M\to I\!\!R$  a transnormal function.
		Then,
		\begin{enumerate}
			\item[(a)] Integral curves of the vector field  $\verb"grad"  \rho$,
			parameterized by arc length, are  geodesics orthogonal to regular leaves.
			\item[(b)] If $(M,g)$ is a  complete Finsler manifold   such that $[a,b] \subset \rho(M)$ does not have critical values, then, for every $p\in \rho^{-1}(b)$,
			$$ d(\rho^{-1}(a),p)=\mathit{l}(\gamma)=d(\rho^{-1}(a),\rho^{-1}(b))= r=\int_{a}^{b}\frac{d s}{\sqrt{\mathfrak{b}(s)} },$$
			where $\gamma:[0,r]\to M$ is the integral curve of  $\verb"grad"  \rho$ parameterized by arc length joining  $\rho^{-1}(a)$ to $p\in\rho^{-1}(b)$, 
			and $\mathit{l}(\gamma)$ is the length of $\gamma$. 
		\end{enumerate}
	\end{prop}
	Note that for a transnormal function  $\rho:M\to[a,b]$  on the complete Finsler manifold $M$ where $a$ and $b$ are the only critical values of $\rho$,  one can extend the geodesic $\gamma$ to $\rho^{-1}(a)$ and $\rho^{-1}(b)$ and so the results of Proposition \ref{proposition-equidistant} can be extended to the whole manifold $M$. That is, we have the following corollary.
\bc \label{Col_1}
 If $(M,g)$ is a  complete Finsler manifold and $\rho: M\to [a,b]$  a transnormal function  such that $a$ and $b$ are the only critical values of $\rho$, then, for every $ c\in[a,b]$ and every $p\in \rho^{-1}(c)$,
 
$$
d(\rho^{-1}(a),p)=\mathit{l}(\gamma)=d(\rho^{-1}(a),\rho^{-1}(c))= r=\int_{a}^{c}\frac{d s}{\sqrt{\mathfrak{b}(s)} }, 
$$
where 	$\gamma:[0,r]\to M$
	 is the unit speed geodesic which joins $\rho^{-1}(a)$ to $p\in \rho^{-1}(c)$ and coincides with the reparametrization of  integral curve of  $\verb"grad"  \rho$ in $\rho^{-1}\{(a,b)\}$. Moreover, the geodesic $\gamma$ is orthogonal to all the leaves $\rho^{-1}(c)$, $c\in[a,b]$.
\ec	
	
	We call the geodesic $\gamma:[0,r]\to M$, that is the unit speed geodesic whose trace coincides with the integral curve of $\verb"grad"  \rho$, a \textit{horizontal geodesic.}

	%%%%%%%%%%%%%%%%
	%%%%%%%%%%%%%%%%%
	
	\section{A Special Solution  of $\nabla^{H}\nabla^{H}\rho= \phi g$  }
	\label{Sec:special}
	Let $(M,g)$ be an n-dimensional Finsler manifold and
	$\rho:M\rightarrow \R$  a solution of Eq.~\eqref{c-field}. If  $\phi$ is a linear  function of $\rho$
	with constant coefficients,  then we say that $\rho$ is a
	\emph{special solution} of Eq.~\eqref{c-field}. Hence, any
	special solution of Eq.~\eqref{c-field} can be written in the form
	\be
	\label{sc-field}
	\nabla^{H}\nabla^{H}\rho= (- K \rho + B) g,
	\ee
	where $K$ and $B$ are constants.
	The  Eq.~\eqref{sc-field} along any geodesic with arc-length $t$ reduces to the ordinary differential equation
	% see akbar-zadeh paper 'sur les spaces de finsler a
	% courbures sectionelles constantes.' proof of '14' in handbook.
	\be \label{equation}
	\frac{d^{2}\rho}{dt^{2}}=-K \rho + B .
	\ee
	Now for the special case $K=C^2>0$ and $B=0$, we have
	\be
	\label{sc-field1} \frac{d^{2}\rho}{dt^{2}}+ C^2\rho=0.
	\ee
	By a
	suitable choice of the arc-length $t$, a solution of Eq.\eqref{sc-field1} is given
	by
	\be \label{rho1}
	\rho (t)= A \cos (Ct),
	\ee and its first derivative is
	\be
	\rho' (t)=-AC \sin (Ct).
	\ee
	So, we can see at a glance  that Eq.~\eqref{rho1}  has   two critical points
	corresponding to $t=0$ and $t=\frac{\pi}{C}$ on $M$ which are repeated  periodically.
	Hence, if $\rho$ is a non-trivial solution of
	Eq.~(\ref{sc-field1}), then it can be written in the following form
	\be \label{rho}
	\rho(t)=\frac{-1}{C}\cos (Ct),\ \quad( A=\frac{-1}{C}).
	\ee
	Taking  Eq.~\eqref{meter} into account,
	%in an adapted coordinates $(u^1=t, u^\alpha,  \
	%\alpha=2,3,...,n )$,
	the  metric form of $M$ becomes
	\be\label{metric form}
	ds^2= dt^2+\big(\sin (Ct)\big)^2\overline{ds}^2,
	\ee
	where $\overline{ds}^2$ is the metric form of a $t$-level of $\rho$
	given by $\overline{ds}^2=f_{\gamma\beta}du^{\gamma}du^{\beta}$. This is the polar form of a Finsler
	metric on a standard sphere of radius $\frac{1}{C}$, see \cite{shen2013differential}.
	
	%%%%%%%%%%%%%%%%%
	%%%%%%%%%%%%%%%%%%
	%%%%%%%%%%%%%%%%%%
	\section{Finsler Manifolds of Positive Constant Flag Curvature} \label{Otava}
	Let  $(x,y)$ be the line element of $TM$ and $P(y,X)\subset T_{x}M$  a 2-plane
	generated by the vectors $y$ and $X$ in
	$T_{x}M $. Then the \emph{flag curvature} $K(x,y,X)$ with respect
	to  the plane $P(y,X)$ at a point $x\in M$ is defined by
	$$
	K(x,y,X):=\frac{g(R(X,y)\,y,X)}{g(X,X)g(y,y)-g(X,y)^{2}},
	$$
	where $R(X,y)\,y$ is the $h$-curvature tensor of Cartan connection. If
	$K$ is independent of $X$, then $(M,g)$ is called \emph{space of}
	\emph{scalar curvature}. If $K$ has no dependence on $x$ or $y$,
	then the Finsler manifold is said to be of \emph{constant (flag) curvature},
	see  for instance \cite{zadeh2006initiation}.  It can be easily verified that the
	components of the $h$-curvature tensor of Cartan connection in the
	adapted coordinate system are given by
	\be \label{2.21}
	\begin{array}{lll}
		R^\alpha_{\ 1 \gamma 1}= - R^\alpha_{\ \gamma 1
			1}=(\frac{\rho'''}{\rho'})\delta_{\gamma}^{\alpha},\\\\
		R^1_{\ 1 \gamma \beta}= - R^ 1_{\ \gamma 1 \beta}=-\rho' \rho'''
		f_{\gamma\beta},\\\\
		R^\alpha_{\ \delta \gamma \beta}=\overline{R}^\alpha_{\ \delta
			\gamma \beta}-(\rho'')^2(f_{\gamma\beta}\delta_{\delta}^{\alpha}-
		f_{\delta\beta}\delta_{\gamma}^{\alpha}),
	\end{array}
	\ee
	where $\overline{R}^\alpha_{\ \delta
		\gamma \beta}$ are components of $h$-curvature tensor related to the metric form $\overline{ds}^2$ on a $t$-level of $\rho$, see \cite{asanjarani2008classification} for more details.
	\bp \label{constant curvature}
	The n-dimensional complete Finsler manifold $(M,g)$  is of
	constant flag curvature $K=C^2>0$, if and only if,
	there is a  non-trivial solution of $\nabla^{H}\nabla^{H}\rho= (- C^2 \rho + B) g$ on $M$.
	\ep
	\bpf
	A Finsler manifold $(M,g)$  is of constant flag
	curvature $K$ if and only if the components of the $h$-curvature tensor are given by  the following, see \cite{asanjarani2008classification} for more details.
	\be \label{K}
	R^i_{\ hjk}= K (\delta^i_hg_{jk}-\delta^i_jg_{hk}).
	\ee
	Using  Eq.~\eqref{K}, we can easily drive the differential equation
	\be \label{bcs}
	\ddot{A}+KAg=0,
	\ee
	where $A$ is the Cartan torsion tensor, $\dot{A}_{ijk}:=(\nabla^H_sA_{ijk})y^s$ and
	$\ddot{A}_{ijk}:=(\nabla^H_s\nabla^H_tA_{ijk})y^sy^t$, see   Section~1.4 of \cite{bao2012introduction} for more details. 
	
	Assume that  $X,
	Y, Z \in \pi^*TM$ are fixed at $v \in I_xM=\{w \in T_xM, \,g(w,w)=1\}$. Let
	$c:I\!\!R
	\rightarrow M$ be the unit-speed geodesic on $(M,g)$ with
	$\displaystyle{{dc} \over {dt}}(0)=v$ and  $\hat c:={dc \over{dt}}$ be the
	canonical lift of $c$ to $TM_0$. Let $X(t)$, $Y(t)$ and $Z(t)$
	denote the parallel sections along $\hat c$ with $X(0)=X$, $Y(0)=Y$
	and $Z(0)=Z$. Put $A(t)=A(X(t),Y(t),Z(t))$, $\dot A(t)=\dot
	A(X(t),Y(t),Z(t))$ and $\ddot A(t)=\ddot A(X(t),Y(t),Z(t))$. Indeed along geodesics, we have $ \displaystyle \frac{d  A}{dt}=\dot A$, $\displaystyle \frac{d \dot A}{dt}=\ddot A$ and  Eq.~\eqref{bcs} becomes
	\be \label{bcs2}
	\frac{d^2A(t)}{dt^2}+K A(t)=0.
	\ee
	The general solution of this differential equation is $A(t)=A_0 \cos \sqrt{K}t+B_0 \sin \sqrt{K}t$.
	Therefore, Eq.~\eqref{sc-field1} which represents a special case of Eq.~\eqref{sc-field} along  geodesics, has a non-trivial solution on $M$.

	Conversely, let $\rho$ given by Eq.~\eqref{rho} be a solution of Eq.~\eqref{sc-field}  on $M$. Then,
	there is an adapted coordinate system on $M$ for which the components of $h$-curvature are given by~\eqref{2.21}. Hence, first and second  equations  of
	\eqref{2.21} satisfy
	\be\label{third }
	R^i_{\ hjk}=
	\frac{-\rho'''}{\rho'}(\delta^i_hg_{jk}-\delta^i_jg_{hk}).
	\ee
	Differentiate  \eqref{rho} with respect to $t$ and replace  the first and third derivatives
	of $\rho$, we obtain $\frac{-\rho'''}{\rho'}=C^2$. Therefore, the first two equations of~\eqref{2.21} satisfy  Eq.~\eqref{K}.
	
	For the third equation of \eqref{2.21}, we recall that as we see in Section~\ref{Sec:special},  $\rho$ has critical points on $M$. Thus,
	from Lemma~\ref{lem}, the
	$t$-levels of $\rho$ are spaces of  positive constant curvature
	$\displaystyle \rho''^2(0)=C^2$. Therefore,  the third equation of \eqref{2.21} becomes
	$$
	R^\alpha_{\ \delta \gamma \beta}=(C^2 -\rho''^2)(f_{\gamma\beta}\delta_{\delta}^{\alpha}-
	f_{\delta\beta}\delta_{\gamma}^{\alpha}).
	$$
	By substituting  $g_{\alpha\beta}= \rho'^2 f_{\alpha\beta}$ and the first and second derivatives of  $\rho$ in the above equation,  we obtain
	$$R^\alpha_{\ \delta \gamma \beta}=C^2(g_{\gamma\beta}\delta_{\delta}^{\alpha}-
	g_{\delta\beta}\delta_{\gamma}^{\alpha}).$$
	So, all three components of Cartan $h$-curvature tensor satisfy Eq.~\eqref{K}  and the Finsler manifold $(M,g)$ is of constant flag curvature
	$K=C^2$.
		\epf

	Now, we are in a position to prove an extension of  \emph{Obata's} theorem
	to Finsler manifolds.
	\bt\label{isometer2}
	Let   $(M,g)$  be a complete connected Finsler manifold of
	dimension $n\geq2$.  Then, $(M,g)$ is isometric to an n-sphere of radius $\displaystyle \frac{1}{C}$ if and only if there is  a non-trivial solution  of the following equation on $M$:
	
	\be \label{sc-field2}
	\nabla^{H}\nabla^{H}\rho+ C^2 \rho g=0.
	\ee
	
	\et
	\bpf
	Let $(M,g)$ be a Finsler manifold which admits  a non-trivial
	solution of  Eq.~\eqref{sc-field2}. According to Proposition~\ref{constant curvature}, $(M,g)$  is of positive constant
	flag curvature $C^2$. So, as we see in Section~\ref{Sec:special},  the  metric form of
	$(M,g)$ is given by~\eqref{metric form} and so $(M,g)$ is isometric to an n-sphere of radius $\displaystyle \frac{1}{C}$.

	Conversely, if $(M,g)$ is isometric to an n-sphere of radius
	$\frac{1}{C}$, then  the metric form of  $M$ is given by
	$ds^{2}=(dt)^{2}+ \sin ^{2}(Ct)\overline{ds}^2$, where
	$\overline{ds}^2$ is the metric form of a hypersurface of $M$.
	This is the  polar form of a
	Finsler
	metric on an $n$-sphere in $I\!\!R^{n+1}$ with the positive  constant curvature $C^2$, see \cite{shen2013differential}.
	Now by substituting the derivative of $\displaystyle  \rho(t)=-\frac{1}{C} \cos (Ct)$  in the metric form of $M$, we obtain  $ds^{2}=(dt)^{2}+ \rho'^2(t)\overline{ds}^2$. Hence,
	$\rho(t)$ is a non-trivial solution of the second order differential equation~\eqref{sc-field1} or equivalently a non-trivial solution of Eq.~\eqref{sc-field2} along geodesics. 
	\epf

	Now, by considering the number of critical points of $\rho$,
	we have the following result.

	\bc \label{iff}
	Let $(M,g)$ be a complete connected Finsler manifold with dimension $n \geq 2$.
	Then, $(M,g)$ is isometrically
	homeomorphic to an $n$-sphere if and only if 
	$\nabla^{H}\nabla^{H}\rho+ C^2\rho
	g=0$ has a non-trivial solution.
	\ec
	\bpf
	Let $(M,g)$ admit a non-trivial solution of  $\nabla^{H}\nabla^{H}\rho+ C^2\rho
	g=0$,
	then from Theorem~\ref{isometer2} we know that it is isometric to an n-sphere of
	radius $\displaystyle \frac{1}{C}$.
	On the other hand, since $M$ is complete, Proposition~\ref{constant curvature} results in  $(M,g)$ is of positive constant curvature. Therefore, by applying the extension   of  \textit{Meyers's theorem}  to Finsler manifolds, see  \cite{akbar1988espaces}, we can conclude that  $M$ is compact.
	Thus, the function $\rho$ admits its absolute maximum and
	minimum values on $M$. Consequently, $\rho$ has two critical points on $M$ and an extension of
	\textit{Milnor theorem} to Finsler geometry, \cite{lehmann1964generalisation}, implies that $(M,g)$ is homeomorphic to an $n$-sphere.
	
	Conversely, let $(M,g)$ be isometrically homeomorphic to an $n$-sphere of radius $\displaystyle \frac{1}{C}$. Then, 
	Theorem~\ref{isometer2}  implies that $\nabla^{H}\nabla^{H}\rho+ C^2\rho g=0$ has  a non-trivial solution on  $M$. 
	\epf

	Following the Obata's theorem in Riemannian geometry a unit
	sphere is characterized by existence of a solution of the differential
	equation $\nabla\nabla f + f g=0$, where $f$ is a certain
	function on Riemannian manifold $(M,g)$ and $\nabla$ is the Levi-Civita connection associated to the Riemannian metric $g$ \cite{gallot1979equations}. Similarly, Theorem \ref{isometer2} implies 
	that in Finsler geometry a unit sphere can be  characterized by existence
	of a solution of  $\nabla^{H}\nabla^{H}\rho+ \rho g=0$, where $\rho$
	is a certain function on Finsler manifold $(M,g)$ and $\nabla^{H}$ is the Cartan horizontal covariant derivative. In analogy with  Riemannian geometry,
	this leads to a definition for an $n$-sphere in Finsler geometry as
	follows.
	\bd
	A  Finslerian $n$-sphere is a  complete connected Finsler manifold   which admits a non-trivial solution of Eq.~\eqref{sc-field2}.
	\ed
	Equivalently, a {\it Finslerian $n$-sphere}  is
	isometrically homeomorphic to an $n$-sphere endowed with a certain Finsler metric.
	\bt
	Let $(M,g)$ be an $n$-dimensional complete connected Finsler manifold. Then, $(M,g)$ has positive constant flag
	curvature $K=C^2$, if and only if, $(M,g)$ is isometrically homeomorphic to an
	$n$-sphere of radius $\displaystyle \frac{1}{C}$ endowed with a certain Finsler metric.
	\et
	\bpf
	Let $(M,g)$ be of positive constant flag curvature $C^2$. As a consequence of
	Proposition~\ref{constant curvature}, there is  a non-trivial solution of Eq.~\eqref{sc-field1} on $M$.
	Thus, by means of Corollary~\ref{iff} it is isometrically homeomorphic to an
	$n$-sphere of radius $\displaystyle \frac{1}{C}$ equipped with the certain Finsler metric form $ds^{2}=(dt)^{2}+ \sin ^{2}(Ct)\overline{ds}^2$.
	
	Conversely, let $(M,g)$ be a  Finsler manifold which is isometrically homeomorphic to an $n$-sphere of radius $\displaystyle \frac{1}{C}>0$. Then, Corollary~\ref{iff} implies that  $M$
	admits a non-trivial solution
	of Eq.~\eqref{sc-field1}. So, from Proposition~\ref{constant
		curvature},  $(M,g)$ is of positive constant flag curvature $C^2$. 
	\epf
	
	\section{Finsler Transnormal Functions}\label{trans}
	Throughout this section, we assume that $\rho:M\to I\!\!R$ is a non-null Finsler transnormal function (see Section \ref{transs}) on  the complete Finsler manifold $(M,g)$. First, we show that there exists an adapted coordinate system on any Finsler manifold that admits a transnormal function. 
	
	\bl \label{trans-eq}
		Let  $\rho: M\to [a,b]$ be  a Finsler transnormal function on  a complete Finsler manifold $(M,g)$ with $g(\verb"grad"  \rho,\verb"grad"  \rho)=\mathfrak{b}\circ \rho$ and no critical values in $(a,b)$. Then, there exists an adapted coordinate system $(u^1=t,u^2,...,u^n)$ on $M$, where $t$ is   parameter of the reparametrization  of integral curve of $\verb"grad"  \rho$. 
	\el
	\bpf
	From Collorally~\ref{Col_1}, one deduces that the reparametrization of integral curve of $\verb"grad" \rho$ is a geodesic of $(M,g)$ and it is orthogonal to every $\rho^{-1}(c)$, for $c\in [a,b]$. Furthermore, from the same corollary, all of these geodesics start from $\rho^{-1}(a)$ and meet  $\rho^{-1}(c)$ at the same time $\int_{a}^{c}\frac{d s}{\sqrt{\mathfrak{b}(s)} }$. Therefore, inspired by Section $1$ of \cite{asanjarani2008classification}, one can consider an adapted coordinate system on the set of regular points of $M$, $M^o:=\rho^{-1}(a,b)$. In other words, there exists a local coordinate system $(u^1=t,u^2,...,u^n)$ on $M^o$ such that $t$ is the parameter of the unit speed geodesic whose trajectory coincides with the integral curve of $\verb"grad" \rho$. In this coordinate system all points belonging to each level set of $\rho$ map into the same value. So, the value of $\rho$ has no dependency on $u^i$, $i=2,...,n$, and just depends on $t$.  Therefore, for every $p\in M^o$,  $\rho(p)=\rho(\gamma(t_1))$, where $\gamma$ is the reparameterization of  integral curve of $\verb"grad"  \rho$ and $t_1$ is the time when $\gamma$ passes through $p$. Since we can extend each horizontal geodesic to   the singular level sets $\rho^{-1}(a)$ and $\rho^{-1}(b)$, while it preserves its properties,  we confirm the existence of the local coordinate system on the whole manifold $M$.  
	\epf
\begin{prop}\label{trans-eq1}	
	Let  $\rho: M\to [a,b]$ be  a Finsler transnormal function on  a complete Finsler manifold $(M,g)$ with $g(\verb"grad"  \rho,\verb"grad"  \rho)=\mathfrak{b}\circ \rho$ and no critical values in $(a,b)$. Then, in adapted coordinate system:	
		\begin{itemize}
			\item [(a)] Level sets of $\rho$ are defined by $t = constant $, where $t\in[0,\int_{a}^{b}\frac{d s}{\sqrt{\mathfrak{b}(s)} }]$, 
			
			\item [(b)] In $\rho^{-1}\{(a,b)\}$,  $\rho$ satisfies the following equation
			\be\label{eq:rho-double}
			\frac{d^2\rho}{dt^2}= \frac12\mathfrak{b}^\prime(\rho),
			\ee
			
			\item [(c)] The Finsler metric form of $M$ is given by $ds^2=(dt)^2+{\rho^\prime}^2 f_{ij}du^i du^j$,  $i,j=2,3, \dots, n$,  where $f_{ij}$ is a Finsler metric tensor on a regular level set of $\rho$ and ${\rho^\prime}^2 f_{ij}$ is the induced metric tensor on this level set.
		\end{itemize}
		%	where $\phi$ is a function of $\rho$.
	\end{prop}
	\bpf	
In an adapted coordinate system on $M$, all points belonging to any level set $\rho^{-1}(c)$, for $c\in [a,b]$, only depend on  $t$. Moreover,   from Corollary~ \ref{Col_1}, all  horizontal geodesics from $\rho^{-1}(a)$ to $\rho^{-1}(c)$ reach to $\rho^{-1}(c)$   at the same time $t=\int_{a}^{c}\frac{d s}{\sqrt{\mathfrak{b}(s)} }$. So we have the proof of  $(a)$  at hand. 

	To prove $(b)$, note that in an adapted coordinate system  on $M$, $\gamma$ is a unitary geodesic whose velocity vector coincides with the positive multiplication of $\verb"grad"  \rho$, i.e. $\gamma'=\frac{\mathbf{grad} \rho}{\sqrt{(\mathbf{grad} \rho, \mathbf{grad} \rho)}}$. So, 
	%	$\frac{t}{F(\verb"grad"  \rho)}$
	\begin{equation*}
		\frac{d^2\rho(\gamma)}{dt^2}= \frac{d}{dt} d\rho\gamma^\prime= \frac{d}{dt} d\rho\frac{\mathbf{grad} \rho}{\sqrt{g(\mathbf{grad}  \rho, \mathbf{grad} \rho)}}=  \frac{d}{dt}\sqrt{\mathfrak{b}(\rho)}= \frac{{\mathfrak{b}^\prime(\rho)}d\rho\gamma^\prime}{2\sqrt{\mathfrak{b}(\rho)}}      =\frac{\mathfrak{b}^\prime(\rho)}{2}, 
	\end{equation*} 
where the third equality comes from Eqs.~\eqref{dis-tra}  and \ref{grad}. \\

	To prove item $(c)$, once we put the adapted coordinate system on $M$,  the metric $g$ can be written as (\ref{metrg}) and  the Finsler metric form is given by \eqref{meter}. Therefore, nothing else is left to be proved.
	\epf

	As a consequence of Lemma~ \ref{trans-eq} and Proposition~\ref{trans-eq1}, one can say that given a complete Finsler manifold $(M,g)$ and a transnormal function $\rho:M\to [a,b]$ on it, $\rho$ is a solution of Eq.~(\ref{para}), at least in the regular part. Therefore, by using results in  \cite{alexandrino2019finsler} and  \cite{asanjarani2008classification}  one can establish several interesting results for Finsler transnormal functions.
	
	\begin{lem}\label{sph.cent}
		Let $(M,g)$ be a complete Finsler manifold and $\rho:M\to [a,b]$ a transnormal function on it.		If $\rho$ has only one critical point $o:=\rho^{-1}(a)$, then each regular level set $\rho^{-1}(c)$  for $c \in (a,b)$ is a hypersphere of radius $r_c=\int_{a}^{c}\frac{d s}{\sqrt{\mathfrak{b}(s)} }$ with center $o$ and constant sectional curvature  $\frac12\mathfrak{b}^\prime(a)  $.  Also, the   Finsler metric form of this level set is $\overline{ds}^2=f_{ij}du^i du^j$,  $i,j=2,3, \dots, n$, where $f_{ij}$ is given by Eq.~\eqref{meter}.
	\end{lem}
	
	\bpf
	From Corollary \ref{Col_1}, given any regular value $c\in(a,b)$, the horizontal geodesics (extensions of integral curves of $\verb"grad"  \rho$) are the geodesics that minimize the distance from $o$ to $\rho^{-1}(c)$. In fact, these geodesics start from $o$ and reach orthogonally to $\rho^{-1}(c)$ at the same time. Therefore, all the points belonging to $\rho^{-1}(c)$ have the same distance $r_c=\int_{a}^{c}\frac{d s}{\sqrt{\mathfrak{b}(s)} }$ from $o$. That means $\rho^{-1}(c)=\{q\in M: \ d(o,q)=r_c\}=S_{r_c}(o)$ in which $S_{r_c}(o)$ is the Finsler sphere of radius ${r_c}$ and center $o$. Also, according to  Lemma~ \ref{trans-eq}, we can consider an adapted coordinate system on $M$. Hence, $\rho$ satisfies Eq.~\eqref{para} which is  equivalent to Eq.~\eqref{c-field} in an adapted coordinate system. Therefore,  from Lemma~\ref{lem}, each $\rho^{-1}(c)$ with metric $\overline{ds^2}=f_{ij}du^i du^j$ is of positive constant sectional curvature $\rho^{''}(0)$. Finally, from item $(b)$ of Proposition \ref{trans-eq1}, $\rho^{''}(0)=\frac12\mathfrak{b}^\prime(a)$.

	\epf
	
	Now, as special case of Proposition \ref{con}, we have the following   classification result on the Finsler manifolds which admit a transnormal function.
	
	\begin{prop}
		Let  $(M,g)$  be a connected complete Finsler manifold of dimension  $n\geq 2$ and $\rho:M\to I\!\!R$ a transnormal function on it. Then, 
		\begin{itemize}
			\item [(a)] If $\rho$ has no critical points, $M$ is conformal to a direct product $J\times M$ of an open interval $J$ of the real line and an $(n-1)$-dimensional complete Finsler manifold $\overline{M}$.
			\item	[(b)] If $\rho$ has one critical point, $M$ is conformal to an $n$-dimensional Euclidean space.
			\item [(c)] If $\rho$ has two critical points, $M$ is conformal to an $n$-dimensional unit sphere in an Euclidean space.
		\end{itemize}
	\end{prop}
	
	For compact Finsler manifolds which admit a transnormal function we have the following theorem.
	\bt\label{sph}
	Let $(M,g)$ be a simply connected and compact Finsler manifold of dimension  $n> 2$ and $\rho:M\to [a,b] $ a transnormal function on it such that $\rho$ has no critical values in $(a,b)$. Then $M$ is homeomorphic to an $n$-sphere.
	\et
	\bpf
	According to Lemma \ref{trans-eq}, there is an adapted coordinate system on $M$ such that  $\rho$ satisfies Eq.~\eqref{para} with $\phi=\frac{\mathfrak{b}^\prime}{2}$. So, $\rho(t)$ has at most two critical points, see Section~3 of \cite{asanjarani2008classification} for more details. Also, since $M$ is compact, $\rho(t)$ takes its maximum and minimum on $M$. Consequently, $\rho$ has exactly two critical points that might be repeated periodically. These critical points  are corresponding to $a$ and $b$. 
	According to the fact that Eq.~(\ref{c-field})  is equivalent to  Eq.~(\ref{para})  in the adapted coordinate system, the transnormal function $\rho$ is a solution of  Eq.~(\ref{c-field}) and the rest of proof is a direct result of Lemma  \ref{sph2}. 
	\epf

	\bp
		Let $(M,g)$ be a complete connected Finsler manifold with dimension $n \geq 2$ and $\rho:M\to [a,b]$ a  transnormal function with $g(\verb"grad"  \rho, \verb"grad"  \rho)=-C^2\rho^2+d$, where $C$ and $d$ are constant positive numbers.	Then, $(M,g)$ is isometrically
		homeomorphic to an $n$-sphere of radius $\frac{1}{C}$.
	\ep
	
	\bpf
	From Proposition~ \ref{trans-eq1},	
	there exists an adapted coordinate system in which
	\be
	\frac{d^2\rho}{dt^2}= \frac12\mathfrak{b}^\prime(\rho)=-C^2\rho.
	\ee 
	So,	the proof is a direct result of Corollary \ref{iff}, where the equation $\nabla^{H}\nabla^{H}\rho+ C^2\rho
	g=0$ reduces to  $\frac{d^2\rho}{dt^2}+ C^2\rho=0$  in the adapted coordinate system.
	\epf

	\subsection{Example}

	To see a simple example illustrating some results of Finsler transormal functions, consider  a calm pond of water that we throw some small piece of stone into it at some time slot. Assume a two-dimensional Euclidean coordinate system on the surface of the pond where the origin is  the point where the stone entered into the water.  The only force  perturbing the water surface is the wind $W(x,y)=\frac13(y,-x)$ blowing across the pond. We want to find the equation of water waves  at each time and also the path equation of water particles (molecules). First, we present the mathematical model of the problem. Assume the open disk $D=\{(x,y)\in \RR^2| x^2+y^2<T\}$, where $T$ is big enough such that $D$ covers the pond. The associated metric to this problem is a special case of Finsler metric which is called \textit{Randers metric }  and is given by  
	$$
	F(y)=\sqrt{\frac{h^2(y,W)+\lambda h(y,y)}{\lambda^2}}-\frac{h(y,W)}{\lambda},
	$$ 
	where   $h$ is the canonical Euclidean metric and $\lambda=1-h(W,W)$, see \cite{shen2001lectures}.   Now we consider the function $\rho:D\too\RR$ defined by $\rho(x,y)=x^2+y^2$. This is not difficult to show that $g(\verb"grad"  \rho ,\verb"grad"  \rho)=2\rho$, where $g$ is the metric with components $(g_{ij})=\left({1 \over 2}
	\left[ \frac{\partial^{2}}{\partial y^{i}\partial y^{j}} F^2
	\right]\right)$. Hence $\rho$ is a transnormal function. Also, it is easy  to show that, for some $\tau< T$,  $\rho^{-1}(\tau)$ coincides with the location of some water wave and therefore the locations of wave are given by preimages of $\rho$, see Section $4.2$ of  \cite{dehkordi2019huygens} for the details. From Lemma  \ref{sph.cent}, each regular level set of $\rho$, that is the location of the  water wave at each time $t=x^2+y^2$, is a circle of radius $r_t=\int_0^t\frac{ds}{\sqrt{2s}}={\sqrt{2t}}$ with center $0$.
	
\begin{figure}[h]
		\centering
		\includegraphics[scale=0.65]{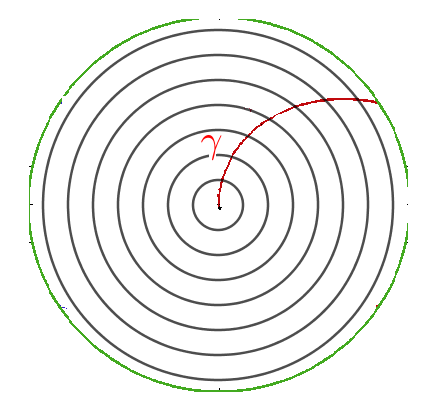} 
		\label{disc.pro}
		\caption{The waves of water and the path of a molecule of water.}
	\end{figure}

	Figure  \ref{disc.pro} illustrates the geodesic $	\gamma(t)=\frac{\sqrt{2}t}{2}(\cos \frac{t}{3}-\sin \frac{t}{3},\sin \frac{t}{3}+\cos \frac{t}{3})$  which is  the track of a molecule of water from time $0$ to time $T$; and also the path of an integral curve of $\verb"grad"  \rho$.  The figure also shows some  $t$-levels of $\rho$, that is the location of waves at different time slots. 
	
	\renewcommand{\baselinestretch}{.5}
	\bibliography{biblio}
	\bibliographystyle{acm}

\end{document}